%% file: PH_Mai.tex
\documentclass[11pt]{article}
\usepackage[T1]{fontenc}
\usepackage{times}
\usepackage{rawfonts}
\usepackage{latexsym}
\usepackage{amsmath}
\usepackage{amsthm}
\usepackage{amsfonts}
\usepackage{amssymb}
\usepackage{enumerate}
\usepackage{xspace}
\usepackage{verbatim}
\usepackage{color}
\newtheorem{theorem}{Theorem}
\newtheorem{lemma}[theorem]{Lemma}
\newtheorem{corollary}[theorem]{Corollary}

\theoremstyle{definition}
\newtheorem{definition}[theorem]{Definition}
\newtheorem{assumption}[theorem]{Assumption}
\newtheorem{remark}[theorem]{Remark}
\newtheorem{example}[theorem]{Example}

\title{Convergence rates of a \vilavmet for nonlinear monotone ill-posed equations in Hilbert spaces}
\author{{\sc Robert Plato}\footnote{Department of Mathematics, University of Siegen,
Walter-Flex-Str.~3, 57068 Siegen, Germany.} $\quad$ and $\quad$ {\sc Bernd Hofmann}%
\footnote{Faculty of Mathematics, Chemnitz University of Technology, 09107 Chemnitz, Germany. 
\newline \hspace*{0.5cm} Research supported by the German Research Foundation (DFG) under grant HO~1454/12-1.}
}

\numberwithin{equation}{section}
\numberwithin{theorem}{section}
\input macros.txt

\makeatletter
\def\thm@space@setup{%
  \thm@preskip=5pt \thm@postskip=5pt
}
\makeatother

\begin{document}
\date{\vspace{-5ex}}
\maketitle
\newcounter{enumcount}
\renewcommand{\theenumcount}{(\alph{enumcount})}
\bibliographystyle{plain}
\begin{abstract}
We consider perturbed nonlinear ill-posed equations in Hilbert spaces, with operators  that are monotone on a given closed convex subset.
A simple stable approach is Lavrentiev regularization, but
existence of solutions of the regularized equation on the given subset can be guaranteed only under additional assumptions that are not satisfied in some applications. 

\lavmet of the related variational inequality seems to be a reasonable alternative then. For the latter approach, in this paper we present new error estimates for suitable a priori parameter choices, if the considered operator is cocoercive and if in addition the solution admits an adjoint source representation.
Some numerical experiments are included.
\end{abstract}
%
%
\section{Introduction}
\label{intro}
In this paper we consider nonlinear equations of the form
\begin{align}
\F u = \fst,
\label{eq:maineq}
\end{align}
where $ \F: \ix \supset \DF \to \ix $ is a nonlinear operator
in a \bhsl{real separable Hilbert space} $ \ix $
with inner product $ \skp{\cdot}{\cdot}: \ix \times \ix \to \reza $,
and $ \fst \in \R(F) = F(\DF) $.
It is assumed that equation \refeq{maineq} is \bhsl{ill-posed} in one of the concepts considered in
\cite{Hofmann_Plato[18]}, \ie it is unstable solvable at
$ \fst $ or locally ill-posed at each solution of \refeq{maineq};
see also \cite{Bot_Hofmann[16]}.
If not specified otherwise, throughout the paper we restrict the considerations to the following class of operators.
\begin{definition}
\label{th:mononote}
The operator $ \F: \ix \supset \DF \to \ix $
is called \emph{monotone} on a set $ \Mset \subset \DF $ if
\begin{align}
\skp{\F\myu-\F\myv}{\myu- \myv} \ge 0 \foreach \;\myu, \myv \in \Mset.
\label{eq:monotone}
\end{align}
\end{definition}
In the following we assume that equation \refeq{maineq} has a solution $ \ust \in \Mset $.
Moreover, we suppose that the \rhs of \refeq{maineq} is only approximately given as $\fdelta \in \ix$
satisfying
\begin{align}
\norm{ \fst - \fdelta } \le \delta,
\label{eq:noisy_data}
\end{align}
where $ \delta \ge 0 $ is a given noise level.
For the regularization of the considered equation \refeq{maineq} with noisy data as in
\refeq{noisy_data},
\lavmet
\begin{align}
(\F + \para I) \uparadeltab = \fdelta,
\label{eq:lavmet}
\end{align}
may be considered, where $ \para > 0 $ is a regularization parameter.
Solvability of equation \refeq{lavmet} on $ \Mset $ is a critical issue and can only be guaranteed under additional assumptions
on the operator $ \F $ and
the set $ \Mset$, \eg
\begin{myenumerate_indent}
\item \label{it:hemi_M_H} $ \F $ is hemicontinuous and $ \DF = \Mset = \ix $, or
\item \label{it:maxmonot} $ \F $ is maximal monotone on $ \Mset $, or
\item \label{it:ball} $ \F $ is hemicontinuous, $ \Mset$ is a closed ball, centered at a solution of \refeq{maineq} and with sufficiently large radius, and $ \tfrac{\delta}{\para} $ is sufficiently small.
\end{myenumerate_indent}
For \ref{it:maxmonot} we refer e.g.~to Deimling~\cite[Theorem 12.5]{Deimling[85]} and note that
\ref{it:hemi_M_H} is a special case of \ref{it:maxmonot} (cf.~e.g.~Showalter~\cite[p.~39]{Showalter[97]}). The case \ref{it:ball} is considered in Tautenhahn~\cite{Tautenhahn[02]}, with some clarification given by Neubauer~\cite{Neubauer[16]}.

There exist examples, however, where none of these conditions
\ref{it:hemi_M_H} -- \ref{it:ball} on $ \F $ and $ \Mset $ is necessarily satisfied.
For other examples, maximal monotonicity in \ref{it:maxmonot} is hard to verify,
e.g.~for operators on $ \ix = L^2(\Omega) $ with $ \Omega \subset \reza^n $, and $ \Mset \subset \{ f \in \ix \mid f \ge 0 \ \textup{a.e.} \}$.


In such cases, a variational formulation (see formula \refeq{lavmet-vi} below) seems to be a reasonable alternative for \refeq{lavmet}. To prove this fact, is one of
the goals of the present paper.

We conclude this section with some references on the regularizing properties of \refeq{lavmet}:
see, \eg
Alber and Ryazantseva~\cite{Alber_Ryazantseva[06]},
Bo\c{t} and Hofmann~\cite{Bot_Hofmann[16]},
Hofmann, Kaltenbacher and Resmerita~\cite{Hofmann_Kaltenbacher_Resmerita[16]},
Janno~\cite{Janno[00]},
Liu and Nashed~\cite{Liu_Nashed[96]},
Tautenhahn~\cite{Tautenhahn[02]}, as well as
Mahale and Nair~\cite{Mahale_Nair[13]}.
\section{\Vilavmet{} -- Basic notations}
\label{varlavmet}
%
%
We introduce the following assumptions and notations.
\begin{assumption}
\label{th:assump_1}
Let $ \F: \ix \supset \DF \to \ix $ be a \bhsl{\demicont bounded operator}
in the real separable Hilbert space $ \ix $ which is \bhsl{monotone} on a given \bhsl{closed convex subset} $ \Mset \subset \ix $, with $ \Mset \subset \DF $.
In addition, let $ \fst, \fdel \in \ix $ satisfy the noise model \refeq{noisy_data}.
Furthermore, we suppose that the equation $ \F \myu = \fst $ has a solution which belongs to $ \Mset $.
\end{assumption}
Throughout the present paper, we assume that Assumption \ref{th:assump_1} holds.
Recall that the operator $ F $ is, by definition in the sense of Deimling~\cite[Definition 11.2]{Deimling[85]} and
Showalter~\cite[p.~36]{Showalter[97]},
\begin{mylist_indent}
\item
\emph{\demicont{}}, if for each $ \myx \in \DF $ and for each sequence $ (\xn) \subset \DF $ with
$ \xn \to \myx $ as $ n \to \infty $, we have weak convergence
$ \F \xn \rightharpoonup \F \myx $ as $ n \to \infty $,

\item \emph{bounded}, if for each bounded set $ \Nset \subset \DF $, the set $ \F(\Nset) \subset \ix $ is bounded.
\end{mylist_indent}
%
%

Instead of \lavmet \refeq{lavmet},
in what follows we consider the following \vilavmet (\ref{eq:lavmet-vi}).
Let, for $ \para > 0$, $  \upardel \in \Mset $ satisfy
\begin{align}
\skp{\F \upardel + \alpha \upardel -\fdel}{\myv-\upardel} \ge 0
\foreach \myv \in \Mset.
\label{eq:lavmet-vi}
\end{align}
For technical purposes, we use for $ \para > 0 $ the notation
\begin{align}
\upara = \uparzer 
\label{eq:upara}
\end{align}
for the noise-free case $\delta=0$, which means that
the approximation obtained by the \vilavmet has been derived on the basis of exact data $ \fdelta = \fst $.
An approach \refeq{lavmet-vi} can be considered as a variational inequality formulation of \lavmet.
%

%
%
A solution to the \varineq \refeq{lavmet-vi} with the penalized operator always exists on $ \Mset $ and depends stably of $ \fdelta $:
%
\begin{theorem}
\label{th:noise-amplific}
Let Assumption \ref{th:assump_1} be satisfied.
Then for each parameter $\para>0 $,
the \vipert \refeq{lavmet-vi} has a unique solution $ \upardel \in \Mset $. In addition, the following stability estimate is satisfied,
\begin{align}
\norm{\upardel - \upara} \le \mfrac{\delta}{\para},
\label{eq:noise-amplific}
\end{align}
where $ \upara \in \Mset  $ is given by \refeq{upara}.
\end{theorem}
\proof
Consider, for $ \para > 0 $ fixed,
the nonlinear operator $ \Fpara: \ix \supset \DF \to \ix$ which maps as $\cdott \myu \mapsto \F \myu + \para \myu $.
Then we obviously have $ \skp{\Fpara \myu - \Fpara \myv}{\myu-\myv} \ge \para \normqua{\myu - \myv} $ for each $ \myu, \cdott \myv \in \Mset $,
\ie
the nonlinear operator $ \Fpara$ is strongly monotone on the set $ \Mset $.
Existence thus follows, \eg from 
Showalter~\cite[proof of Theorem 2.3 in Chapter II]{Showalter[97]}.
%
The mentioned proof in that reference may be applied, using the notations from there, with $ \mathcal{A} = \Fpara $ and $ v_0 = \ust $, where again $ \ust \in \Mset $ satisfies $ \F \ust = \fst $. Notice that the operator $ \mathcal{A} $ considered 
in \cite{Showalter[97]} is assumed to be pseudo-monotone all over the considered Hilbert space.  However, the proof in that paper can be employed straightforward under the assumptions made in the present paper.

%

We next verify estimate \refeq{noise-amplific}. For notational convenience we introduce the notation
\begin{align*}
\mydiffdel = \upardel - \upara \foreach \para > 0.
\end{align*}
We have
\begin{align*}
\skp{\F \upara + \para \upara - \fst}{ \mydiffdel } \ge 0,
\quad
\skp{\F \upardel + \para \upardel - \fdelta}{ -\mydiffdel } \ge 0.
\end{align*}
Summation of those two inequalities gives
\begin{align*}
0 & \le
\skp{\F \upara + \para \upara - \fst}{ \mydiffdel }
-\skp{\F \upardel + \para \upardel - \fdelta}{ \mydiffdel }
\\
&
=
-\skp{\F \upardel -\F \upara}{ \mydiffdel }
-\para \skp{\mydiffdel}{\mydiffdel}
+ \skp{\fdelta -\fst}{ \mydiffdel }
\le
0
-\para \normqua{\mydiffdel}
+ \delta \norm{\mydiffdel},
\end{align*}
and the statement of the theorem follows by rearranging terms.
\proofend
\begin{remark}
Existence results for \varineqs (either for similar, more general or more specific situations) may also be found in others papers and monographs. See, e.g.~Barbu and Precupanu~\cite[Theorem 2.67 and subsequent remark]{Barbu_Precupanu[12]}
and Kinderlehrer and Stampacchia~\cite[Corollary 1.8 of Chapter III]{Kinderlehrer_Stampacchia[00]}.
In Bakushinsky, Kokurin and Kokurin~\cite[Lemma 6.1.3]{Bakushinsky_Kokurin_Kokurin[18]}
and Br\'{e}zis~\cite[Proposition 31]{Brezis[68]},
a simple proof is given for the special case that the operator $ \F $ satisfies a Lipschitz condition on the monotonicity set $ \Mset $.

Using more refined arguments, it is possible to weaken the assumptions of Theorem \ref{th:noise-amplific} without changing the statement of the theorem:
the condition ``separable'' on the Hilbert space $ \ix $ can be removed in fact, and the assumption 
``\demicont, bounded'' on the operator $ \F $ can be replaced by the weaker property ``hemicontinuous''; \cf Browder~\cite{Browder[65]}
or Br\'{e}zis~\cite[Theorem 24]{Brezis[68]}.
\remarkend
\end{remark}

\bn
Below we consider the overall regularization error $ \upardel - \ust $, where $ \ust \in \Mset $ denotes a classical or generalized solution of the equation $ \F \myu = \fst $,
\cf \refeq{maineq} above or
\refeq{vi-unpert} below.
This overall error can be decomposed into regularization error
$ \upara - \ust $ and noise amplification term $ \upardel - \upara $. The latter term has already been estimated in \refeq{noise-amplific}, and we thus have
\begin{align}
\norm{\upardel - \ust} \le \norm{\upara - \ust} + \mfrac{\delta}{\para}
\;\foreach \;\para > 0.
\label{eq:regularization-error}
\end{align}
%
%
%
%
%
Below we thus may focus on the estimation of the bias norm
$ \norm{\upara - \ust} $.
\section{Convergence of regularized solutions}
\label{convergence}
In this section, we
consider strong convergence of the elements $ \upara $
generated by the \vilavmet 
\refeq{lavmet-vi} 
as $ \para \to 0 $.
We continue to assume that the conditions stated in Assumption \ref{th:assump_1} are satisfied.

As a preparation, we consider the unperturbed, unpenalized version of the \vilavmet \refeq{lavmet-vi}, i.e.~the determination of an
$ \vst \in \Mset $ which satisfies
the \varineq
\begin{align}
\skp{\F \vst -\fst}{\myv-\vst} \ge 0 \foreach \myv \in \Mset.
\label{eq:vi-unpert}
\end{align}
\begin{remark}
\label{th:varineq-comments}
\begin{myenumerate}

\item
We note that any classical solution of \refeq{maineq} obviously satisfies
the \varineq \refeq{vi-unpert}, so the set of solutions of
\refeq{vi-unpert} is by assumption non-empty.

\item
The \varineq \refeq{vi-unpert} is equivalent with
$ \skp{\F \myv -\fst}{\myv-\vst} \ge 0 $ for each $ \myv \in \Mset $,
\cf \eg Showalter~\cite[Corollary~2.4]{Showalter[97]}),
or Browder~\cite[Lemma 1]{Browder[65]}.
This in particular means that the set of solutions satisfying the \varineq \refeq{vi-unpert} is closed and convex, and thus it has a unique element with minimal norm $ \ustst $.

\item
Let the operator $ F $ be strictly monotone on $ \Mset $, \ie
in \refeq{monotone} we may replace ``$ \ge $'' by strict inequality ``$>$'' for each $ \myu, \myv \in \Mset $ with $ \myu \neq \myv $. Then \refeq{vi-unpert} and also
\refeq{maineq} have at most one solution, respectively.

\item
Any element $ \vst \in \Mset $ solves
the \varineq \refeq{vi-unpert} if and only if
the identity
$ \vst = \PM (\vst - \mu \kla{\F \vst-\fst}) $ holds for each $ \mu \ge 0 $, where
$ \PM: \ix \to \ix $ denotes the convex projection onto the set $ \Mset $.
This follows from a standard variational formulation for convex projections, see \eg
Kinderlehrer and Stampacchia~\cite[Theorem 2.3 of Chapter I]{Kinderlehrer_Stampacchia[00]}.
A similar statement holds for the \vilavmet \refeq{lavmet-vi}.
\remarkend
\end{myenumerate}
\end{remark}
\begin{theorem}
\label{th:convergence}
Let Assumption \ref{th:assump_1} be satisfied.
We have $ \upara  \to \ustst $ as $ \para \to 0 $, where $ \ustst \in \Mset $ denotes the minimum norm solution of the \varineq \refeq{vi-unpert}.
\end{theorem}
\proof
This easily follows, e.g., by a compilation of the steps considered in 
the proof of Theorem 3 in Ryazantseva~\cite{Ryazantseva[76]}.
\proofend
\begin{remark}
Convergence of the \vilavmet is in fact the subject of many research papers and monographs, see \eg Alber and Ryazantseva~\cite[Theorem 4.1.1]{Alber_Ryazantseva[06]},
Bakushinsky, Kokurin and Kokurin~\cite[Lemma 6.1.4]{Bakushinsky_Kokurin_Kokurin[18]},
Khan, Tammer and Zalinescu~\cite{Khan_Tammer_Zalinescu[15]},
Liu and Nashed~\cite{Liu_Nashed[98]}, and
Ryazantseva~\cite{Ryazantseva[76]},
and the references therein. 

Quite frequently in the literature, more general situations than in the present paper are considered, \eg perturbation of the considered convex set $ \Mset $ in \refeq{vi-unpert}, or set-valued operators $ \F $ in Banach spaces. On the other hand, the assumptions made in Theorem \ref{th:convergence} are weaker in some aspects. For example, we allow the monotonicity set in \refeq{monotone} to be a nontrivial subset of $\ix $, with a possibly empty interior, and in addition
no Lipschitz continuity of the operator $ \F $ is required
in Theorem \ref{th:convergence}.
\remarkend
\end{remark}

\bn
As an immediate consequence of Theorem \ref{th:convergence} and estimate
\refeq{regularization-error}, we obtain the following result.
\begin{corollary}
Let Assumption \ref{th:assump_1} be satisfied.
For any a priori parameter choice $ \para = \para(\delta) $
with $ \para(\delta) \to 0 $ and
$ \tfrac{\delta}{\para(\delta)} \to 0 $ as $ \delta \to 0 $, we have
%
\begin{align}
\upardeldel \to \ustst
\quad \textup{as } \  \delta \to 0,
\end{align}
where $ \ustst $ is as in Theorem \ref{th:convergence}.
\end{corollary}
\section{Convergence rates for regularized solutions}
\label{convergence_rates}
In this section, we provide convergence rates of $ \upara $ as $ \para \to 0 $ under
adjoint source conditions.
%
We continue to assume that the conditions stated in Assumption \ref{th:assump_1} are satisfied.
In addition, the following class of operators will be of importance, \cf Bauschke and Combettes~\cite[Definition 4.4]{Bauschke_Combettes[85]}.
\begin{definition}
\label{th:cocoercive}
An operator $ \F: \ix \supset \DF \to \ix $ in a Hilbert space $ \ix $ is called \emph{\cocoercive} on  a subset $ \Mset \subset \DF $ if, for some constant
$ \tau > 0 $, we have
\begin{align}
\skp{\F \myu - \F \myv}{\myu-\myv} \ge \tau \normqua{\F \myu - \F  \myv}
\foreach \myu, \cdott \myv \in \Mset.
\label{eq:cocoercive}
\end{align}
\end{definition}
A \cocoercive operator is sometimes called inverse strongly monotone.
For $ \tau > 0 $ fixed, an operator $ \F $ is \cocoercive on $ \Mset $ with constant $   \tau $ if and only if $ I - \mu F $ is nonexpansive for each $ 0 \le \mu \le  2\tau $.
Cocoerciveness obviously implies monotonicity. An example of a \cocoercive operator may be found in Liu and Nashed~\cite[Example 3]{Liu_Nashed[98]}. Another example is given in section~\ref{paraesti} of the present paper.

%
%
Below, frequently we make use of the following Lipschitz condition.
\begin{assumption}
\label{th:assump_2}
Let $ \DF \subset \ix $ be an open subset,
and let $ \F $ be \frechet differentiable on $ \DF $.
In addition, let the following Lipschitz condition be satisfied on a given subset
$ \Mset \subset \DF $,
\begin{align}
\norm{\prim{F}(\myu)-\prim{\F}(\myv)} \le L \norm{\myu-\myv}
\foreach \myu, \myv \in \Mset,
\label{eq:lipschitz-prime}
\end{align}
where $ L \ge 0 $ denotes some finite constant.
\end{assumption}
%
The following proposition provides a useful tool for the verification of cocoerciveness of a nonlinear operator.
\begin{proposition}
\label{th:coco-diffable}
Let Assumptions \ref{th:assump_1} and \ref{th:assump_2}
be satisfied. Let $ \prim{F}(\myu) $ be \cocoercive on $ \ix $, uniformly for
$ \myu \in \Mset $,
%
\ie there exists some constant $ \tau > 0 $ such that
for each $ \myu \in \Mset $
\begin{align}
\skp{\prim{F}(\myu)h}{h} \ge \tau \normqua{\prim{F}(\myu)h} \quad
\forall h \in \ix,
\label{eq:coco-diffable}
\end{align}
holds. Then $ \F $ is \cocoercive on $ \Mset $, with constant $ \tau $.
\end{proposition}
\proof
From uniform cocoerciveness of $ \prim{F} $, we obtain
for any $ \myu \in \Mset $ and $ h \in \ix $ with $ \myu + h \in \Mset $ that
$ \F(\myu+h) - \F(\myu) = \inttxt{0}{1}{ \prim{\F}(\myu+th)h}{dt} $,
and thus
\begin{align*}
& \skp{\F(\myu+h) - \F(\myu)}{h}
= \ints{0}{1}{ \skp{\prim{\F}(\myu+th)h}{h}}{dt}
\ge \tau \ints{0}{1}{ \normqua{\prim{\F}(\myu+th)h}}{dt}
\\
&
\quad
\ge
\tau \klabi{\ints{0}{1}{ \norm{\prim{\F}(\myu+th)h}}{dt}}^2
\ge
\tau \normqua{\ints{0}{1}{ \prim{\F}(\myu+th)h}{dt}}
= \tau
\normqua{\F(\myu+h) - \F(\myu)}.
\proofend
\end{align*}
\begin{remark}
\label{th:monotone-diffable}
\begin{myenumerate}
\item
If $ \prim{F}(\myu) $ is a monotone operator on $ \ix $ for each
$ \myu \in \Mset $, then $ F $ is monotone on $ \Mset $.
This immediately follows from the proof of Proposition \ref{th:coco-diffable}
by considering the case $ \tau = 0 $ there.

\item It is evident from the proof of Proposition \ref{th:coco-diffable} that in \refeq{coco-diffable}, ``$ \forall h \in \ix $'' can be replaced by the weaker condition ``$ \forall h \in \ix $ satisfying $ \myu + t h \in \Mset $ for $ t > 0 $ sufficiently small'', without changing the statement of the proposition.
One can show that this in fact yields an equivalent condition for \cocoerciveness.
\remarkend
\end{myenumerate}
\end{remark}
%
For ill-posed problems, convergence rates can only be obtained under additional conditions on the solution. In this section we assume that
there exists a solution of equation \refeq{maineq} which belongs to $ \Mset $ and satisfies an adjoint source condition, \ie
\begin{align}
\ust \in \Mset, \quad
\F \ust = \fst, \quad
\ust = \Fprime(\ust)^* z, \quad \norm{z} =: \varrho,
\label{eq:adjoint-source-condition}
\end{align}
for some $ z \in \ix $.
This completes the formulation of the basic assumptions needed in this section.

For the proof of the main result of this section, \cf Theorem \ref{th:epara-speed} below,
we need the following lemma.
For any element $ \myu \in \Mset $ consider
\begin{align}
\eparau \defeq \eparau(\myu) = \upara - \myu, \quad
\rpara = \F \upara - \fst, \quad
\epara = \eparau(\ust)
\for \para > 0,
\label{eq:difference-notations}
\end{align}
where $ \upara \in \Mset $ is introduced in \refeq{upara}.
\begin{lemma}
\label{th:vi-lemma}
Let Assumption \ref{th:assump_1} be satisfied.
For any $ \myu \in \Mset $ we have,
with the notations from
\refeq{difference-notations},
\begin{align}
\skp{\rpara}{\eparau} + \para \normqua{\eparau} \le - \para \skp{\myu}{\eparau}
\for \para > 0.
\label{eq:vi-lemma}
\end{align}
%
\end{lemma}
\proof
We consider \refeq{lavmet-vi} with $ \delta = 0 $, which means $ \fdel = \fst $ in fact:
\begin{align*}
\skp{\F \upara - \fst + \para \upara}{\upara -\myu}
=
\skp{\rpara + \para \upara}{\eparau}
=
\skp{\rpara}{\eparau} +  \para \skp{\upara}{\eparau} \le 0.
\end{align*}
%
%
From this we obtain
\begin{align*}
\skp{\rpara}{\eparau} + \para \normqua{\eparau}
=
\skp{\rpara}{\eparau} +  \para \skp{\upara}{\eparau}
-\para \skp{\myu}{\eparau}
\le - \para \skp{\myu}{\eparau},
\end{align*}
which is \refeq{vi-lemma}.
\proofend

\bn
We are now in a position to formulate the main result of this section.
\begin{theorem}
\label{th:epara-speed}
Let Assumptions \ref{th:assump_1} and \ref{th:assump_2}
be fulfilled.
If $ \F $ is \cocoercive on $ \Mset $, and if in addition the adjoint source condition \refeq{adjoint-source-condition} is satisfied with
$ \varrho L < 2 $, then
\begin{align}
\norm{ \upara - \ust } &= \Landauno{\para^{1/2}},
\qquad
\norm{\F \upara - \fst} = \Landauno{\para} \as \para \to 0.
\label{eq:epara-rpara-speed}
\end{align}
\end{theorem}
\proof
We proceed with \refeq{vi-lemma} for $ \myu = \ust $.
From \refeq{adjoint-source-condition} we obtain, with the notations introduced in
\refeq{difference-notations},
\begin{align}
-\skp{\ust}{\epara} =
-\skp{\Fprime(\ust)^* z}{\epara} =
-\skp{z}{\Fprime(\ust) \epara}
\le
\varrho \norm{\Fprime(\ust) \epara}.
\label{eq:epara-speed-a}
\end{align}
For a further estimation of \refeq{epara-speed-a}, we need to consider the first order remainder $ \R = \R_{\ust} $ of a Taylor expansion at $ \ust \in \DF $:
\begin{align*}
\R(\myu) &= \F(\myu)- \F(\ust) - \prim{\F}(\ust)\kla{\myu-\ust},
\quad \myu \in \DF.
\end{align*}
For $ h \in \ix $ such that the line segment from $ \ust $ to $ \ust + h $ belongs to $ \DF $,
we have
$ \R(\ust+h) = \inttxt{0}{1}{ \kla{ \prim{\F}(\ust+th)-\prim{\F}(\ust)}h}{dt} $ and thus
%
$ \norm{\R(\ust+h)} \le \tfrac{L}{2} \normqua{h} $.
%
This gives
$ \Fprime(\ust) \epara = \F(\upara) - \F(\ust) - \R(\upara)
= \rpara - \R(\upara) $
with $ \norm{\R(\upara)} \le \tfrac{L}{2} \normqua{\epara} $.
We are now in a position to proceed with the upper bound in
\refeq{epara-speed-a}:
\begin{align}
\norm{\Fprime(\ust) \epara } \le \norm{\rpara} + \norm{\R(\upara)}
\le \norm{\rpara} + \mfrac{L}{2} \normqua{\epara}.
\label{eq:epara-speed-b}
\end{align}
The estimates \refeq{vi-lemma}
for $ \myu = \ust $ and
\refeq{epara-speed-a} -- \refeq{epara-speed-b}
finally give
\begin{align*}
\skp{\rpara}{\epara} + \para \normqua{\epara}
\le - \para \skp{\ust}{\epara}
\le \varrho \para \norm{\Fprime(\ust) \epara}
\le \varrho \para \klabi{\norm{\rpara} + \mfrac{L}{2} \normqua{\epara}},
\end{align*}
and thus
\begin{align}
\skp{\rpara}{\epara} + \para \klabi{1- \tfrac{\varrho L}{2}} \normqua{\epara}
\le \varrho \para \norm{\rpara}.
\label{eq:epara-speed-c}
\end{align}
This in particular means
$ \skp{\rpara}{\epara} \le \varrho \para \norm{\rpara} $, and
\cocoerciveness, \cf \refeq{cocoercive}, moreover means
$ \skp{\rpara}{\epara} \ge \tau \normqua{\rpara} $.
We thus obtain
\begin{align}
\tau \norm{\rpara} \le \varrho \para,
\label{eq:epara-speed-d}
\end{align}
\ie
$ \norm{\rpara} = \Landauno{\para} $ as $ \para \to 0 $.
From \refeq{epara-speed-c} and \refeq{epara-speed-d} we finally obtain
\begin{align*}
\tau \klabi{1- \tfrac{\varrho L}{2}} \normqua{\epara}
\le \tau \varrho \norm{\rpara}
\le \varrho^2 \para,
\end{align*}
which is the first statement in \refeq{epara-rpara-speed}.
\proofend
\begin{remark}
\begin{myenumerate}
\item
From Theorem \ref{th:epara-speed} and Theorem \ref{th:convergence}
it follows that any $ \ust $ satisfying the conditions in \refeq{adjoint-source-condition} is the minimum norm solution of the \varineq \refeq{vi-unpert}.

\item
Theorem \ref{th:epara-speed} improves results in Liu and Nashed~\cite[Theorem 6]{Liu_Nashed[98]}, where $ \norm{\eparalong} = \Landauno{\para^{1/3}} $ as $ \para \to 0 $ is obtained only (under more general assumptions, however, \eg possible set perturbations).

\item The first error estimate in Theorem \ref{th:epara-speed} remains valid if in
\refeq{adjoint-source-condition}, the identity $ \F \ust = \fst $ is replaced by the weaker assumption that $ \ust \in \Mset $ satisfies the \varineq \refeq{vi-unpert}. In the proof of Theorem \ref{th:epara-speed},
then one only has to make additional use of the fact that
the inequality 
$ \skp{\F \upara - \F \ust }{\eparalong} \le  \skp{\F \upara - \fst }{\eparalong} $ holds.
The second error estimate in Theorem~\ref{th:epara-speed} 
has to be replaced by $ \norm{\F \upara - F \ust} = \Landauno{\para} $ then.

\item
Using some ideas of Tautenhahn~\cite{Tautenhahn[02]} and
Janno~\cite{Janno[00]}, one may obtain convergence rates
for source conditions of the form
$ \ust = \Fprime(\ust) z $, \ie the
adjoint source condition is replaced by a classical one.
This topic, however, goes beyond the scope of the present study and will be considered elsewhere.

\item
For recent results on adjoint source conditions for linear problems, see Plato, Hofmann, and Math\'{e}~\cite{Plato_Hofmann_Mathe[16]}.
\remarkend
\end{myenumerate}
\end{remark}
\begin{corollary}
\label{th:apriori}
Under the conditions of Theorem \ref{th:epara-speed} we have,
for any a priori parameter choice $\paradelta \sim \delta^{2/3}$,
the convergence rate result
\begin{align}
\norm{u_{\paradelta}^\delta - \ust}
= \Landau(\delta^{1/3}) \quad \textup{ as }\;\; \delta \to 0.
\label{eq:apriori}
\end{align}
\end{corollary}
\begin{remark}
\begin{myenumerate}
\item
The rate \refeq{apriori} is identical with rates obtained in
\cite[Theorem 3, Remark 4]{Hofmann_Kaltenbacher_Resmerita[16]} for \lavmet \refeq{lavmet}
with variational source conditions.

\item
The rate of convergence in \refeq{apriori} is higher than those obtained by
Liu and Nashed~\cite{Liu_Nashed[98]},
Thuy~\cite{Thuy[11]}, and
Buong~\cite{Buong[05]} for the \vilavmet under similar source conditions. Note that, on the other hand, the results in those papers are established in a more general framework, respectively, \eg in Banach spaces or allowing set perturbations, and for a~posteriori parameter choice strategies.
%
\remarkend
\end{myenumerate}
\end{remark}
\section{Modified \vilavmet}
\label{lavmet-translate}
Occasionally it may be useful to consider a modified version of the \vilavmet \refeq{lavmet-vi}.
For this purpose let $ \myubar \in \ix $ be fixed.
For $ \para > 0 $ let $ \upardel \in \Mset $ satisfy
\begin{align}
\skp{\F \upardel + \alpha( \upardel - \myubar) -\fdel}{\myv-\upardel} \ge 0
\foreach \myv \in \Mset.
\label{eq:lavmet-vi-translate}
\end{align}
We denote by $ \upara = \uparzer $ the approximation obtained by the
modified \vilavmet \refeq{lavmet-vi-translate} with exact data $ \fdelta = \fst $.
Method \refeq{lavmet-vi-translate} can be considered as variational inequality formulation of the  translated \lavmet $ \F \uparadeltab + \para (\uparadeltab - \myubar) = \fdelta $.
The results of sections \ref{varlavmet}--\ref{convergence_rates} can be easily applied to the modified \vilavmet by considering translation: replace the operator $ \F $ and the monotonicity set $ \Mset $
there by
\begin{align*}
\Ftil:
\ix \supset  - \myubar + \DF \to \ix, \ v \mapsto  \F(\myubar + v),
\quad
\Mtil = - \myubar + \Mset,
\end{align*}
respectively. We briefly formulate the relevant results under 
the general assumption that the conditions stated in Assumption \ref{th:assump_1} are satisfied.
\begin{myenumerate_indent}
\item
The modified \vilavmet \refeq{lavmet-vi-translate} has a unique solution $ \upardel \in \Mset  $ which satisfies $ \norm{\upardel - \upara} \le \tfrac{\delta}{\para} $ for each $ \para > 0 $.
\item
We have $ \upara  \to \ustst $ as $ \para \to 0 $, where $ \ustst \in \Mset $ denotes the solution of the \varineq \refeq{vi-unpert} having minimal distance to $ \myubar $.
In addition, for any a priori parameter choice
$ \para = \paradelta $
with $ \paradelta \to 0 $ and
$ \tfrac{\delta}{\paradelta} \to 0 $ as $ \delta \to 0 $, we have
$ \upardeldel \to \ustst $ as $ \delta \to 0 $.

\item
If Assumption \ref{th:assump_2} is fulfilled and $ \F $ is \cocoercive on $ \Mset $, and if in addition the adjoint source condition
\begin{align}
\ust \in \Mset, \quad
\F \ust = \fst, \quad
\ust - \myubar = \Fprime(\ust)^* z, \quad
\varrho \defeq \norm{z},
\label{eq:adjoint-source-cond-mod}
\end{align}
is satisfied with some $ z \in \ix $ and
$ \varrho L < 2 $, then
\begin{align*}
\norm{\upara - \ust} = \Landauno{\para^{1/2}} \ \textup{ as } \para \to 0,
\qquad
\norm{\upardeldel - \ust} = \Landau(\delta^{1/3}) \ \textup{ as } \delta \to 0,
\end{align*}
for any a priori parameter choice $\paradelta \sim~\delta^{2/3}$.
\end{myenumerate_indent}
%
An appropriate choice of $ \myubar $ guarantees that $ \ust - \myubar $ belongs to the range of $ \Fprime(\ust)^* $, which typically requires, besides sufficient smoothness, that appropriate conditions on a subset of the boundary of the domain of definition $ \DF $ are satisfied.
%
%
%
%
%
%
%
\section{An example, and numerical illustrations}
\subsection{A parameter estimation problem}
\label{paraesti}
We consider the estimation of the coefficient $ u \in L^2(0,1) $ in the following        initial value problem:
\begin{align*}
\prim{f} + u f = 0 \; \text{\ a.e.~on } [0,1], \quad f(0) = -\czer < 0,
\end{align*}
where $ f \in H^1(0,1) $;
cf.~Groetsch~\cite{Groetsch[93]}, Hofmann~\cite{Hofmann[99]}, or
Tautenhahn~\cite{Tautenhahn[02]}.
The initial value $ -\czer $ with $ \czer > 0 $ is assumed to be known exactly.
This problem can be written as
$ F u = f $, with
\begin{align}
(F u)(\myt) \defeq -\czer e^{-U(\myt)},
\quad U(\myt) = \ints{0}{\myt}{ u(\mys) } { d\mys }, \quad 0 \le \myt \le 1.
\label{eq:fex}
\end{align}
The operator $ F : L^2(0,1) \to L^2(0,1) $ is bounded and \frechet differentiable on
$ L^2(0,1) $, with \frechet derivative
\begin{align}
[\prim{\F}(u)h](t) = - (\F u)(t) H(t) \for
h \in L^2(0,1), \
H(\myt) = \ints{0}{\myt}{ h(\mys) } { d\mys }, \ 0 \le \myt \le 1.
\label{eq:fex-prime}
\end{align}
Let 
\begin{align}
\Dsetc = \inset{ u \in L^2(0,1) \mid u \ge \mykapa\; \text{\ a.e.~on}\;[0,1] }, 
\label{eq:Muzer}
\end{align}
where
$ \mykapa \in \reza $.
\begin{proposition}
\label{th:Fu-ex-frechet-mono-coco}
The operator $ \F $ in \refeq{fex} is monotone on $ \Dsetc[0] $.
For any $ \mykapa > 0 $, it is \cocoercive
on $ \Dsetc $, with constant $ \tau = \tfrac{\mykapa}{2\czer} $.
\end{proposition}
\proof
We shall make use of Proposition \ref{th:coco-diffable} and Remark \ref{th:monotone-diffable}. Let $ \myu \in L^2(0,1), \fex \defeq -Fu $, and $ h \in L^2(0,1) $.
From \refeq{fex-prime} it follows that
\begin{align*}
\skp{\prim{F}(u)h}{h}
& =
\ints{0}{1}{ (\fex H) \prim{H} } { d\myt }
=
\fex H^2 \vert_0^1
- \ints{0}{1}{ \prim{(\fex H)} H}{ d\myt }
\ge - \ints{0}{1}{ \prim{\fex} H^2 }{ d\myt }
-\skp{\prim{F}(u)h}{h},
\end{align*}
and thus
\begin{align}
2 \skp{\prim{F}(u)h}{h}
\ge - \ints{0}{1}{ \prim{\fex} H^2 }{ d\myt }
= \ints{0}{1}{ \fex \myu H^2 }{ d\myt },
\label{eq:Fu-ex-monotone-a}
\end{align}
where the properties $ \fex(1) H^2(1) \ge 0 $ and $ H(0) = 0 $ have been used.
Estimate \refeq{Fu-ex-monotone-a} implies for each $ \myu \in \Dsetc[0] $ that
$ \skp{\prim{F}(u)h}{h} \ge 0 $ for each $ h \in \ix $,
and the monotonicity statement for $ \F $ immediately follows from
Remark \ref{th:monotone-diffable}.
%

Now let $ \mykapa > 0 $ be fixed. For any $ \myu \in \Dsetc $
we proceed with \refeq{Fu-ex-monotone-a}:
%
%
%
\begin{align*}
2 \skp{\prim{F}(u)h}{h}
\ge
\mykapa \ints{0}{1}{ \fex H^2 }{ d\myt }
\ge
\tfrac{\mykapa}{\czer} \ints{0}{1}{ (\fex H)^2 }{ d\myt }
= \tfrac{\mykapa}{\czer} \ints{0}{1}{ (\prim{F}(u)h)^2 }{ d\myt }
= \tfrac{\mykapa}{\czer} \normqua{\prim{F}(u)h},
\end{align*}
where the estimate $ \fex \le \czer $ has been applied.
The cocoerciveness statement for $ \F $ now follows from
Proposition \ref{th:coco-diffable}.
\proofend
\subsection{Numerical experiments}
The theoretical results are finally illustrated by some numerical experiments
for the operator $ F : L^2(0,1) \to L^2(0,1) $
considered in \refeq{fex}, with $ \czer = 1 $ there.
We give a few preparatory notes on the numerical tests first.
\begin{mylist}
\item
In each of our numerical experiments we choose a convex closed subset $ \Mset = \Dsetc $ of the form
\refeq{Muzer}
%
%
with some lower bound $ \mykapa > 0 $.
The setting \refeq{Muzer} guarantees
\cocoerciveness (cf.~Proposition~\ref{th:Fu-ex-frechet-mono-coco}),
and Lipschitz continuity \refeq{lipschitz-prime} of the operator
$ \Fprime $ on $ \Mset $ holds with $ L = \czer = 1 $.
We consider some $ \ust \in H^1(0,1) $ with $ \ust \in \Mset $, and then
the adjoint source condition \refeq{adjoint-source-cond-mod} is satisfied for $ \myubar \equiv \ust(1) $.
The solution $ \ust $ and the set $ \Mset $ are always chosen in such a way that the condition $ \varrho L < 2 $ is satisfied,
\cf \refeq{adjoint-source-cond-mod} and the subsequent conclusion there.

\item
We consider the a priori parameter choice $ \paradelta = \delta^{2/3} $, for different values of $ \delta $.

\item
The modified \vipert \refeq{lavmet-vi-translate}
is approximately solved by using a fixed point iteration for the corresponding
fixed point equation
\begin{align*}
\uparadelta = \PM (\uparadelta- \mu \kla{\F \uparadelta + \para (\uparadelta - \myubar) -\fdelta}),
\with \para = \paradelta,
\end{align*}
and the initial guess is the function $ \myubar $. 
Notice that the underlying fixed point operator is contractive, with contraction constant
$ 1-\mu \para $,
provided that the step size satisfies $ 0 < \mu < 2\tau = \mykapa $, \cf the remarks following Definition \ref{th:cocoercive}, and Proposition \ref{th:coco-diffable}.
In addition, the regularization parameter must satisfy
$ 0 < \para \le \tfrac{1}{\mu} - \tfrac{1}{\mykapa} $.
In our numerical experiments we always choose $ \mu = \tfrac{\mykapa}{2} $.
%

Iteration is stopped if the norm difference of two consecutive iterates satisfies, for the first time, an estimate of the form $ \le c \delta $, with some constant $ c > 0 $. This stopping criterion ensures that the resulting approximation $ \upardeldeltil \in \Mset $ satisfies
$ \norm{\upardeldeltil - \upardeldel} = \Landau(\delta^{1/3}) $,
which is of sufficient accuracy.

\item
The problem is discretized
using a backward rectangular rule for the integrals,
and replacing each considered (continuous) function $ \psi:[0,1] \to \reza $ by
$ (\psi(nh))_{n=0,\ldots, N} $,
 with step size $ h = \tfrac{1}{N} $ for $ N = 200 $. This leads to a fully discretized nonlinear problem in $ \reza^{N+1} $.

\item
In the numerical experiments we consider perturbations of the form $ f_n^\delta = f(nh) + \Delta_n, \ n = 0,1,\ldots, N $, with uniformly distributed random values $ \Delta_n $ satisfying $ \modul{\Delta_n} \le \delta $.
\end{mylist}
\begin{example}
\label{th:example1}
We first consider the equation $ \F u = \fst $,
with \rhs
\begin{align*}
\fst(\myt) = -\exp(-\tfrac{\myalpha}{2} \myt^2 -\mybeta \myt) \for 0 \le \myt \le 1,
\end{align*}
with $ \myalpha = \mybeta = \tfrac{1}{2} $. The exact solution is then given
by
\begin{align*}
\ust(\myt) = \myalpha \myt + \mybeta \for 0 \le \myt \le 1.
\end{align*}
We may consider the set $ \Mset = \Dsetc $ in \refeq{Muzer} with $ \mykapa = \mybeta $.
%
The numerical results are given in Table \ref{tab:num1}.
\begin{table}[h]
\hfill
\begin{tabular}{|| c | c |@{\hspace{5mm} } c | c ||}
\hline
\hline
$ \delta $
& $ 100 \myast \delta/\norm{f} $
& $ \ \norm{\upardeldel - \ust} $
& $ \ \norm{\upardeldel - \ust} \ / \delta^{1/3} \ $
\\ \hline \hline
 $1.0 \myast 10^{-2}$ & $1.33 \myast 10^{0}$ & $9.87 \myast 10^{-2}$ & $0.46$ \\
 $5.0 \myast 10^{-3}$ & $6.66 \myast 10^{-1}$ & $8.23 \myast 10^{-2}$ & $0.48$ \\
 $2.5 \myast 10^{-3}$ & $3.33 \myast 10^{-1}$ & $6.72 \myast 10^{-2}$ & $0.50$ \\
 $1.2 \myast 10^{-3}$ & $1.67 \myast 10^{-1}$ & $5.42 \myast 10^{-2}$ & $0.50$ \\
 $6.2 \myast 10^{-4}$ & $8.33 \myast 10^{-2}$ & $4.17 \myast 10^{-2}$ & $0.49$ \\
 $3.1 \myast 10^{-4}$ & $4.16 \myast 10^{-2}$ & $3.26 \myast 10^{-2}$ & $0.48$ \\
 $1.6 \myast 10^{-4}$ & $2.08 \myast 10^{-2}$ & $3.26 \myast 10^{-2}$ & $0.61$ \\
 $7.8 \myast 10^{-5}$ & $1.04 \myast 10^{-2}$ & $2.72 \myast 10^{-2}$ & $0.64$ \\
 $3.9 \myast 10^{-5}$ & $5.21 \myast 10^{-3}$ & $2.53 \myast 10^{-2}$ & $0.75$ \\
\hline
\hline
\end{tabular}
\hfill \mbox{}
\caption{Numerical results for Example \ref{th:example1}}
\label{tab:num1}
\end{table}
\remarkend
\end{example}

\begin{example}
\label{th:example2}
We next consider the equation $ \F u = \fst $
with \rhs
\begin{align*}
\fst(\myt) = -\exp(\tfrac{\myalpha}{\pi} (\cos \pi \myt -1) - \mybeta \myt)
\for 0 \le \myt \le 1,
\end{align*}
with $ \myalpha = \tfrac{1}{4}, \
\mybeta = \tfrac{1}{3} $. The exact solution is then given
by
\begin{align*}
\ust(\myt) = \myalpha \sin \pi \myt + \mybeta \for 0 \le \myt \le 1.
\end{align*}
We may consider the set $ \Mset = \Dsetc $ in \refeq{Muzer} with 
$ \mykapa = \mybeta $.
%
The numerical results are shown in Table \ref{tab:num2}.
\begin{table}[h]
\hfill
\begin{tabular}{|| c | c |@{\hspace{5mm} } c | c ||}
\hline
\hline
$ \delta $
& $ 100 \myast \delta/\norm{f} $
& $ \ \norm{\upardeldel - \ust} $
& $ \ \norm{\upardeldel - \ust} \ / \delta^{1/3} \ $
\\ \hline \hline
 $1.0 \myast 10^{-2}$ & $1.25 \myast 10^{0}$ & $7.00 \myast 10^{-2}$ & $0.32$ \\
 $5.0 \myast 10^{-3}$ & $6.25 \myast 10^{-1}$ & $4.66 \myast 10^{-2}$ & $0.27$ \\
 $2.5 \myast 10^{-3}$ & $3.12 \myast 10^{-1}$ & $3.87 \myast 10^{-2}$ & $0.29$ \\
 $1.2 \myast 10^{-3}$ & $1.56 \myast 10^{-1}$ & $3.01 \myast 10^{-2}$ & $0.28$ \\
 $6.2 \myast 10^{-4}$ & $7.81 \myast 10^{-2}$ & $2.22 \myast 10^{-2}$ & $0.26$ \\
 $3.1 \myast 10^{-4}$ & $3.90 \myast 10^{-2}$ & $1.60 \myast 10^{-2}$ & $0.24$ \\
 $1.6 \myast 10^{-4}$ & $1.95 \myast 10^{-2}$ & $1.08 \myast 10^{-2}$ & $0.20$ \\
 $7.8 \myast 10^{-5}$ & $9.76 \myast 10^{-3}$ & $7.54 \myast 10^{-3}$ & $0.18$ \\
 $3.9 \myast 10^{-5}$ & $4.88 \myast 10^{-3}$ & $4.70 \myast 10^{-3}$ & $0.14$ \\
\hline
\hline
\end{tabular}
\hfill \mbox{}
\caption{Numerical results for Example \ref{th:example2}}
\label{tab:num2}
\end{table}
\remarkend
\end{example}

\bibliography{PH}

\end{document}

%% file: macros.txt
\newtheorem{proposition}[theorem]{Proposition}
\newcommand{\para}{\alpha}
\newcommand{\paradelta}{\para(\delta)}
\newcommand{\ix}{\mathcal{H}}

\newcommand{\upara}[1][\para]{u_{#1}}
\newcommand{\uparadelta}[1][\para]{u_{#1}^\delta}
\newcommand{\uparadeltab}[1][\para]{u}

\newcommand{\fdelta}{f^\delta}

\newcommand{\reza}{ \mathbb{R} }
\newcommand{\Landau}{\mathcal{O}}
\newcommand{\refeq}[1]{(\ref{eq:#1})}

\newcommand{\F}{F}

\newcommand{\normqua}[1]{\norm{#1}^2}

\newcommand{\skp}[2]{\langle #1,  #2 \rangle}
\newcommand{\norm}[1]{\Vert \hspace{0.4mm} #1 \hspace{0.4mm} \Vert}

\theoremstyle{nonumberplain}

\def\endproof{\vbox{\hrule height0.4pt\hbox{%
   \vrule height0.75ex width0.5pt\hskip0.8ex
   \vrule width0.5pt}\hrule height0.4pt
   }}

\newenvironment{myenumerate}{%
\begin{list}{(\alph{enumcount})}
{\setcounter{enumcount}{1}\usecounter{enumcount}
\setlength{\topsep}{1mm}
\setlength{\itemsep}{0mm}
\setlength{\labelwidth}{0mm}
\setlength{\labelsep}{1mm}
\setlength{\itemindent}{1mm}
\setlength{\leftmargin}{0mm}
}}{\end{list}}

\newcommand{\DF}{\Dset(\F)}

\newcommand{\Dset}{\mathcal{D}}
\newcommand{\Mset}{\mathcal{M}}
\newcommand{\Nset}{\mathcal{N}}
\newcommand{\myu}{u}
\newcommand{\myv}{v}
\newcommand{\myx}{v}
\newcommand{\ust}{u_*}
\newcommand{\vst}{\ust}
\newcommand{\ustst}{u_{**}}

\newcommand{\fst}{f_*}
\newcommand{\fdel}{f^\delta}
\newcommand{\upardel}{u_\para^\delta}
\newcommand{\upardeldel}{u_{\paradelta}^\delta}
\newcommand{\upardeldeltil}{\widetilde{u}_{\paradelta}^\delta}
\newcommand{\uparzer}{u_\para^0}
\newcommand{\foreach}{\quad \text{for each} \ \ }
\newcommand{\with}{\quad \text{with} \ \ }
\newcommand{\lavmet}{Lavrentiev regularization\xspace}
\newcommand{\vipert}{penalized \varineq{}\xspace}
\newcommand{\vilavmet}{penalized \varineq{} method\xspace}
\newcommand{\Vilavmet}{Penalized \varineq{} method\xspace}

\newcommand{\tinybullet}{{\tiny \raisebox{0.6mm}{$ \bullet $}}}
\newenvironment{mylist}{%
\begin{list}{\tinybullet}
{\setlength{\topsep}{0.2cm}
\setlength{\itemsep}{0mm}
\setlength{\labelwidth}{0mm}
\setlength{\labelsep}{3mm}
\setlength{\itemindent}{3mm}
\setlength{\leftmargin}{0mm}
}}{\end{list}}

\newenvironment{mylist_indent}{%
\begin{list}{\tinybullet}
{\setlength{\topsep}{0.2cm}
\setlength{\itemsep}{-1mm}
\setlength{\labelwidth}{2mm}
\setlength{\labelsep}{3mm}
\setlength{\itemindent}{0mm}
\setlength{\leftmargin}{5mm}
}}{\end{list}}

\newcommand{\demicont}{demicontinuous\xspace}
\newcommand{\demiconty}{demicontinuity\xspace}
\newcommand{\varineq}{variational inequality\xspace}

\newcommand{\varineqs}{variational inequalities\xspace}

\newcommand{\klabi}[1]{\big(\cdottsm #1 \cdottsm\big)}
\newcommand{\cdott}{\hspace{0.4mm}}
\newcommand{\cdottsm}{\hspace{0.3mm}}
\newcommand{\kla}[1]{(#1)}
\newcommand{\mfrac}[2]{\dfrac{\mbox{\footnotesize \raisebox{-0.5mm}{$#1$}}}%
{\mbox{\footnotesize \raisebox{0.8mm}{$#2$}}}}
\newcommand{\proofend}{\qquad \endproof}
\newcommand{\frechet}{Fr\'{e}chet\xspace}
\newcommand{\Fprime}{\prim{F}}
\newcommand{\prim}[1]{#1^{\prime}}
\newcommand{\cocoercive}{cocoercive\xspace}
\newcommand{\cocoerciveness}{\cocoercive{}ness\xspace}
\newcommand{\epara}{e_\para}
\newcommand{\eparalong}{u_\para -\ust}
\newcommand{\eparau}{\Delta_{\para}}
\newcommand{\Landauno}[1]{\Landau\kla{#1}}
\newcommand{\as}{\quad \text{as} \ \ }
\newcommand{\rpara}{r_\para}
\newcommand{\Fpara}{\F_\para}
\newcommand{\R}{\mathcal{R}}
\newcommand{\cf}{cf.\mbox{}\xspace}
\newcommand{\ie}{i.e.,\xspace}
\newcommand{\eg}{e.g.,\xspace}
\newcommand{\remarkend}{\quad \ensuremath{\vartriangle}}
\newcommand{\bn}{\bigskip \noindent}
\newcommand{\defeq}{:=}
\newcommand{\for}{\quad \text{for} \ \ }

\newcommand{\xn}[1][n]{\myv_{#1}}
\newcommand{\ints}[4]{\mathop{\raisebox{-0.1mm}{\mbox{\Large$ \textstyle \int $}}}\nolimits_{\hspace{-1mm}#1}^{#2} \hspace{-0.2mm} #3 \, #4}
\newcommand{\inttxt}[4]{\int_{#1}^{#2}\hspace{-0.2mm} #3 \, #4}
\newcommand{\rhs}{right-hand side\xspace}
\newcommand{\paraa}{\para}
\newcommand{\parab}{\beta}
\newcommand{\uparaa}{u_{\paraa}}
\newcommand{\uparab}{u_{\parab}}
\newenvironment{myenumerate_indent}{%
\begin{list}{(\alph{enumcount})}
{\setcounter{enumcount}{1}\usecounter{enumcount}
\setlength{\topsep}{1mm}
\setlength{\itemsep}{-0.4mm}
\setlength{\labelwidth}{5mm}
\setlength{\labelsep}{2mm}
\setlength{\itemindent}{-0mm}
\setlength{\leftmargin}{7mm}
}}{\end{list}}
\newcommand{\PM}{\mathcal{P}_{\Mset}}
\newcommand{\mydiffdel}{\chi_\para^\delta}
\newcommand{\myubar}{\overline{u}}

\newcommand{\fex}{g}
\newcommand{\inset}[1]{\{ \, #1 \, \}}
\newcommand{\Dsetc}[1][\mykapa]{\mathcal{M}_{#1}}

\newcommand{\Ftil}{\widetilde{\F}}
\newcommand{\Mtil}{\widetilde{\Mset}}
\newcommand{\myast}{\cdot}
\newcommand{\mykapa}{\theta}
\newcommand{\modul}[1]{\vert \hspace{0.3mm} #1 \hspace{0.3mm} \vert}
\newcommand{\myt}{t}
\newcommand{\mys}{s}
\newcommand{\myalpha}{a}
\newcommand{\mybeta}{b}
\newcommand{\bhsl}[1]{#1}
\newcommand{\czer}{c_0}

%% file: PH_Mai.bbl
\begin{thebibliography}{10}

\bibitem{Alber_Ryazantseva[06]}
Y.~Alber and I.~Ryazantseva.
\newblock {\em Nonlinear Ill-posed Problems of Monotone Type}.
\newblock Springer-Verlag, Berlin, Heidelberg, New York, 1st edition, 2006.

\bibitem{Bakushinsky_Kokurin_Kokurin[18]}
A.~B. Bakushinsky, M.~M. Kokurin, and M.~Yu Kokurin.
\newblock {\em Regularization Algorithms for Ill-Posed Problems}.
\newblock de Gruyter, Berlin, 1st edition, 2018.

\bibitem{Barbu_Precupanu[12]}
V.~Barbu and T.~Precupanu.
\newblock {\em Convexity and Optimization in {B}anach spaces}.
\newblock Springer, Dordrecht, Heidelberg, London, New York, 4th edition, 2012.

\bibitem{Bauschke_Combettes[85]}
H.~H. Bauschke and P.~L. Combettes.
\newblock {\em Convex Analysis and Monotone Operator Theory}.
\newblock Springer-Verlag, Berlin, Heidelberg, New York, 2nd edition, 2017.

\bibitem{Bot_Hofmann[16]}
R.~I. Bo\c{t} and B.~Hofmann.
\newblock Conditional stability versus ill-posedness for operator equations
  with monotone operators in {H}ilbert space.
\newblock {\em Inverse Problems}, 32(12):125003 (23pp), 2016.

\bibitem{Brezis[68]}
H.~Br\'{e}zis.
\newblock \'{E}quations et in\'{e}quations non lin\'{e}ares dans les espaces
  vectoriels en dualit\'{e}.
\newblock {\em Annales de l'institut Fourier}, 18(1):115--175, 1968.

\bibitem{Browder[65]}
F.~E. Browder.
\newblock Nonlinear monotone operators and convex sets in {B}anach spaces.
\newblock {\em Bull. Amer. Math. Soc.}, 71(5):780--785, 1965.

\bibitem{Buong[05]}
N.~Buong.
\newblock Convergence rates in regularization for ill-posed variational
  inequalities.
\newblock {\em CUBO}, 7(3):87--94, 2005.

\bibitem{Deimling[85]}
K.~Deimling.
\newblock {\em Nonlinear Functional Analysis}.
\newblock Springer-Verlag, Berlin, Heidelberg, New York, 1st edition, 1985.

\bibitem{Groetsch[93]}
C.~W. Groetsch.
\newblock {\em Inverse Problems in the Mathematical Sciences}.
\newblock Vieweg, Braunschweig, Wiesbaden, 1993.

\bibitem{Hofmann[99]}
B.~Hofmann.
\newblock {\em Mathematik Inverser Probleme}.
\newblock Teubner, Stuttgart, Leipzig, 1 edition, 1999.

\bibitem{Hofmann_Kaltenbacher_Resmerita[16]}
B.~Hofmann, B.~Kaltenbacher, and E.~Resmerita.
\newblock Lavrentiev's regularization method in {H}ilbert spaces revisited.
\newblock {\em Inverse Problems and Imaging}, 10(3):741--764, 2016.

\bibitem{Hofmann_Plato[18]}
B.~Hofmann and R.~Plato.
\newblock On ill-posedness concepts, stable solvability and saturation.
\newblock {\em Journal of Inverse and Ill-Posed Problems}, 26(2):287--297,
  2018.

\bibitem{Janno[00]}
J.~Janno.
\newblock Lavr'entev regularization of ill-posed problems containing nonlinear
  near-to- monotone operators with application to autoconvolution equation.
\newblock {\em Inverse Problems}, 16(2):333--348, 2000.

\bibitem{Khan_Tammer_Zalinescu[15]}
A.~A. Khan, C.~Tammer, and C.~Zalinescu.
\newblock Regularization of quasi-variational inequalities.
\newblock {\em Optimization}, 64(8):1703--1724, 2015.

\bibitem{Kinderlehrer_Stampacchia[00]}
D.~Kinderlehrer and G.~Stampacchia.
\newblock {\em An Introduction to Variational Inequalities and Their
  Applications}.
\newblock SIAM, Philadelphia, 1st, reprint edition, 2000.

\bibitem{Liu_Nashed[96]}
F.~Liu and M.~Z. Nashed.
\newblock Convergence of regularized solutions of nonlinear ill-posed problems
  with monotone operators.
\newblock In {\em Partial Differential Equations and Applications. Lecture
  Notes in Pure and Applied Mathematics Vol. 177}, pages 353--361, New York,
  1996. Marcel Dekker.

\bibitem{Liu_Nashed[98]}
F.~Liu and M.~Z. Nashed.
\newblock Regularization of nonlinear ill-posed variational inequalities and
  convergence rates.
\newblock {\em Set-Valued Analysis}, 6:313--344, 1998.

\bibitem{Mahale_Nair[13]}
P.~Mahale and T.~Nair.
\newblock Lavrentiev regularization of nonlinear ill-posed equations under
  general source conditions.
\newblock {\em Journal of Nonlinear Analysis and Optimization}, 4(2):193--204,
  2013.

\bibitem{Neubauer[16]}
A.~Neubauer.
\newblock Private communication, 2016.

\bibitem{Plato_Hofmann_Mathe[16]}
R.~Plato, B.~Hofmann, and P.~Math\'{e}.
\newblock Optimal rates for {L}avrentiev regularization with adjoint source
  conditions.
\newblock {\em Math. Comp.}, 87:785--801, 2018.

\bibitem{Ryazantseva[76]}
I.~P. Ryazantseva.
\newblock Regularization of non-linear equations with monotonic discontinuous
  operators.
\newblock {\em U.S.S.R. Comput. Math. Math. Phys.}, 16:228--232, 1976.

\bibitem{Showalter[97]}
R.~E. Showalter.
\newblock {\em Monotone Operators in Banach Space and Partial Differential
  Equations}.
\newblock AMS, Providence, Rhode Island, 1st edition, 1997.

\bibitem{Tautenhahn[02]}
U.~Tautenhahn.
\newblock On the method of {L}avrentiev regularization for nonlinear ill-posed
  problems.
\newblock {\em Inverse Problems}, 18:191--207, 2002.

\bibitem{Thuy[11]}
N.~T.~T. Thuy.
\newblock Regularization of ill-posed mixed variational inequalities with
  non-monotone perturbations.
\newblock {\em J. Inequalities Appl.}, 2011:25, 2011.

\end{thebibliography}
